
\documentstyle{amsppt}
\topmatter
\author S.~Giller
\footnote[*]{supported by KBN grant 2P30221706p02 \hfill{ }},
 C.~Gonera$^*$, P.~Kosi\' nski$^*$, P.~Ma\' slanka$^*$
\endauthor
\title Differential calculus on deformed
$\operatorname{E(2)}$ group
\endtitle
\address
S.~Giller, C.~Gonera, P.~Kosi\' nski
\newline Department of Field Theory
\newline University of \L \' od\' z
\newline Pomorska 149/153
\newline 90-236 \L \' od\' z, Poland
\newline
\newline
P.~Ma\' slanka
\newline Institute of Mathematics
\newline University of \L \' od\' z
\newline Banacha 22
\newline 90-238 \L \' od\' z, Poland
\endaddress
\abstract
Fourdimensional bicovariant differential $*$-calculus
on quantum $\operatorname{E(2)}$
\linebreak
group is constructed.
The relevant  Lie algebra is obtained and covariant
differential calculus on quantum plane is found.
\endabstract
\endtopmatter
\NoRunningHeads
\document
\TagsOnRight
\define\di{\operatorname{d}}
\define\ti{\widetilde}
\head I.~Introduction \endhead

Much attention has been paid recently to the special
deformation of Poincar\' e algebra (group) called
$\varkappa$-deformed Poincar\' e algebra \cite{1}
(group \cite{2}). Many particular results were
obtained which gave more deep insight into
their structure. One of the most important open
problems is to find and classify bicovariant
differential $*$-calculi. Partial results were
obtained in Ref.~\cite{3} where some differential
calculi on deformed Minkowski space were considered using
bicrossproduct technique developed in \cite{4}.

In the present paper, following the approach developed
by Woronowicz school \cite{5}, \cite{6} we construct
fourdimensional bicovariant $*$-calculus for the
deformed $\operatorname{E(2)}$ group. This group
differs trivially from twodimensional $\varkappa$-Poincar\' e
group so our results apply, {\it mutatis mutandis},
to this case also. In the subsequent paper we shall
consider the bicovariant differential calculi
for fourdimensional case.

The paper is organized as follows. In section II we describe
the bicovariant differential $*$-calculus obtained, according
to Woronowicz \cite{5}, by an appropriate choice of right
ideal $\Cal R \subset \ker \varepsilon$. In section III the
generalized Lie algebra is found and compared with the
deformed (euclidean) $\varkappa$-Poincar\' e algebra. It is
shown that Woronowicz functionals are expressible in terms of
$\varkappa$-Poincar\' e algebra; in particular, additional,
fourth generator, is related to the Casimir operator of
$\varkappa$-Poincar\' e. Finally, in section IV the covariant
differential calculus on $\varkappa$-plane is constructed and
compared with that given in Ref.~\cite{3}.

To conclude the Introduction let us remind the definition of
deformed $\operatorname{E(2)}$ group \cite{7}, \cite{8}.
It is generated by four elements: $a$, $a^*$, $v_+$,
$v_-=v_+^*$ subject to the following commutation rules~:
$$\aligned
[a,v_-]=\frac{i}{\varkappa}(I-a),\qquad
& [a^*,v_-]=\frac{i}{\varkappa} \left(a^*-{a^*}^2 \right), \\
[a,v_+]=\frac{i}{\varkappa}\left( a-a^2 \right) ,\qquad
& [a^*,v_+]=\frac{i}{\varkappa} \left(I-a^* \right), \\
[v_+,v_-]=\frac{i}{\varkappa}(v_--v_+),\qquad
& a^*a=aa^*=I.
\endaligned \tag1$$
The coproduct, antipode and counit read, respectively~:
$$\aligned
\triangle a=a\otimes a,\qquad
& \triangle a^*=a^*\otimes a^*, \\
\triangle v_+=a\otimes v_+ +v_+\otimes I,\qquad
& \triangle v_-=a^*\otimes v_- +v_-\otimes I, \\
S(a)=a^*,\qquad & S(a^*)=a, \\
S(v_+)=-a^* v_+, \qquad & S(v_-)=-av_-, \\
\varepsilon (a)=\varepsilon (a^*)=1, \qquad
& \varepsilon (v_-)=\varepsilon (v_+)=0.
\endaligned \tag2$$

\head II.~The bicovariant $*$-calculus \endhead

We follow closely the Woronowicz school approach \cite{5}, \cite{6}.
The ideal $\Cal R\subset \ker \varepsilon$ is chosen to be
generated by the following nine elements~:
$$\aligned
(a-I)^2, \qquad & (a^* -I)^2,\\
(a-I)v_+, \qquad & (a^* -I)v_+,\\
(a-I)v_-, \qquad & (a^* -I)v_-,\\
(a-I)(a^*-I)= & 2I-(a+a^*), \\
v_+^2-\frac{i}{\varkappa}v_+,
\qquad &
v_-^2-\frac{i}{\varkappa}v_- .
\endaligned \tag3$$
Actually, some of the generators are redundant but we include them
for simplicity.

We easily check, that $\Cal R$ is $\operatorname{ad}$-invariant;
also $x\in \Cal R$ implies $S(x)^* \in \Cal R$. Therefore,
$\Cal R$ determines bicovariant $*$-calculus \cite{5}. This calculus
is fourdimensional: $\ker \varepsilon /\Cal R$ is spanned by
$a-a^*$, $v_+$, $v_-$, $\dsize v_+ v_- +\frac{i}{\varkappa} v_+$.

Using standard technique we arrive at the following commutation rules
between the elements of algebra and differentials~:
$$\aligned
a\di a=\di a a, \quad a^* \di a=\di a a^*,\quad
& a\di a^*=\di a^* a, \quad a^*\di a^*=\di a^* a^*, \\
v_+\di a=\di a v_+ +\frac{i}{\varkappa}(a-I)\di a,\qquad
& v_+\di a^*=\di a^* v_+ +\frac{i}{\varkappa}(I-a)\di a^*, \\
v_-\di a=\di a v_- +\frac{i}{\varkappa}(a^*-I)\di a,\qquad
& v_-\di a^*=\di a^* v_- +\frac{i}{\varkappa}(a^*-I)\di a^*, \\
a\di v_+=\di v_+ a -\frac{i}{\varkappa} a\di a, \qquad
& a\di v_-=\di v_- a -\frac{i}{\varkappa} a^*\di a, \\
a^*\di v_+=\di v_+ a^* -\frac{i}{\varkappa} a\di a^*, \qquad
& a^*\di v_-=\di v_- a^* -\frac{i}{\varkappa} a^*\di a^*, \\
v_+ \di v_+=\di v_+ v_+ -\frac{i}{\varkappa} \di v_+ , \qquad
& v_- \di v_-=\di v_- v_- -\frac{i}{\varkappa} \di v_- , \\
v_+ \di v_-=-\di v_+ v_- -\frac{i}{\varkappa}\di v_+
+\di & \left( v_+ v_- +\frac{i}{\varkappa}v_+ \right), \\
v_- \di v_+=-\di v_- v_+ -\frac{i}{\varkappa}\di v_-
+\di & \left( v_+ v_- +\frac{i}{\varkappa}v_+ \right).
\endaligned \tag4$$

The left invariant forms read
$$\aligned
& \omega_0=\frac{1}{2} (a^* \di a-a\di a^*)=a^*\di a=a\di a^*, \\
& \omega_+=a^* \di v_+, \\
& \omega_- =a\di v_-, \\
& {\widetilde \omega}_0=\di \left( v_+ v_-
+\frac{i}{\varkappa} v_+ \right) -v_+ \di v_- -v_-\di v_+ .
\endaligned \tag5$$

As a next step we calculate the right action of $\operatorname{E(2)}$
on $\omega$'s~:
$$\aligned
& _{\Gamma} \triangle (\omega_0) =\omega_0
\otimes I, \\
& _{\Gamma} \triangle ({\widetilde \omega}_0) ={\widetilde \omega}_0
\otimes I, \\
& _{\Gamma} \triangle (\omega_+) =\omega_+
\otimes a^* + \omega_0\otimes a^* v_+, \\
& _{\Gamma} \triangle (\omega_-) =\omega_-
\otimes a - \omega_0\otimes a v_- .
\endaligned \tag6$$
Acordingly, the right-invariant forms read
$$\aligned
& \eta_0=\omega_0, \\
& {\widetilde \eta}_0={\widetilde \omega}_0, \\
& \eta_+=a\omega_+ +\left( \frac{i}{\varkappa} I-v_+ \right) \omega_0, \\
& \eta_-=a^*\omega_- +\left( -\frac{i}{\varkappa} I+v_- \right) \omega_0.
\endaligned \tag7$$
Relations (4)--(7) lead, via Woronowicz construction \cite{5}, to the
following external algebra~:
$$\aligned
& \omega_0 \wedge \omega_0=0, \\
& {\widetilde \omega}_0 \wedge {\widetilde \omega}_0 =0, \\
& \omega_0 \wedge \omega_{\pm} =-\omega_{\pm} \wedge \omega_0, \\
& \omega_0 \wedge {\widetilde \omega}_0
=-{\widetilde \omega}_0 \wedge \omega_0, \\
& \omega_{\pm} \wedge {\widetilde \omega}_0
=-{\widetilde \omega}_0 \wedge \omega_{\pm}
\pm \frac{2}{\varkappa^2} \omega_0 \wedge \omega_{\pm}, \\
& \omega_{\pm} \wedge \omega_{\pm}=\mp \frac{i}{\varkappa}
\omega_0 \wedge \omega_{\pm}, \\
& \omega_- \wedge \omega_+=-\omega_+ \wedge \omega_-
+\frac{i}{\varkappa} \omega_0 \wedge (\omega_- -\omega_+).
\endaligned \tag8$$
Also, using equations (4), (5) and (8) one easily gets
$$\aligned
& \di \omega_0=0, \\
& \di {\widetilde \omega}_0=0, \\
& \di \omega_{\pm}=\mp \omega_0 \wedge \omega_{\pm}.
\endaligned \tag9$$
Finally, the $*$-structure is given by the following relations~:
$$\aligned
& \omega_0^*=-\omega_0, \\
& {\widetilde \omega}_0^*=
- {\widetilde \omega}_0, \\
& \omega_{\pm}^*=\omega_{\mp} \pm \frac{i}{\varkappa} \omega_0.
\endaligned \tag10$$

\head III.~ The Lie algebra of deformed
$\operatorname{E(2)}$ group \endhead

It has been shown recently \cite{7}, \cite{8}, \cite{9} that
the deformed Lie algebra $\frak e\operatorname{(2)}$, dual to
$\operatorname{E(2)}$, has the following structure~:
$$\aligned
& [P_1,P_2]=0, \\
& [J,P_1]=i P_2, \\
& [J,P_2]=-i\varkappa \sinh \left( \frac{P_1}{\varkappa} \right) ,
\endaligned \tag11a$$
$$J^*=J, \quad P_1^*=P_1, \quad P_2^*=P_2, \tag11b$$
$$\aligned
& \triangle P_1=I\otimes P_1 +P_1\otimes I, \\
& \triangle P_2=e^{-\frac{P_1}{2\varkappa}}\otimes P_2
+P_2\otimes e^{\frac{P_1}{2\varkappa}}, \\
& \triangle J=e^{-\frac{P_1}{2\varkappa}}\otimes J
+J\otimes e^{\frac{P_1}{2\varkappa}}.
\endaligned \tag11c$$
On the other hand, Woronowicz theory \cite{5} provides us
with the general construction of Lie algebra once bicovariant
calculus is given. The resulting Lie algebra is quadratic so it must
differ from one given by equations (11). Therefore the question
arises how the both algebras are related to each other.
To address this problem we first construct the Woronowicz
algebra for $\operatorname{E(2)}$. To this end we introduce
counterparts of left-invariant vector fields by the formula
$$\di x=(\chi_0 *x)\omega_0 +({\widetilde \chi}_0 * x)
{\widetilde \omega}_0 +(\chi_+ * x)\omega_+
+(\chi_- *x)\omega_- ,\tag12$$
where $x\in \operatorname{E(2)}$ and
$\chi *x=(\operatorname{id} \otimes \chi ) \triangle (x)$.

Using $\di (\di x)=0$ and equations (8), (9) and (12) we arrive
at the following algebra~:
$$\aligned
& [\chi_0,{\widetilde \chi}_0]=0, \\
& [\chi_{\pm},{\widetilde \chi}_0]=0, \\
& [\chi_+,\chi_-]=0, \\
& [\chi_0,\chi_+]=\chi_+ +\frac{i}{\varkappa}\chi_+^2
+\frac{i}{\varkappa} \chi_- \chi_+
-\frac{2}{\varkappa^2}\chi_+ {\widetilde \chi}_0 ,\\
& [\chi_0,\chi_-]=-\chi_- -\frac{i}{\varkappa}\chi_-^2
-\frac{i}{\varkappa} \chi_- \chi_+
+\frac{2}{\varkappa^2}\chi_- {\widetilde \chi}_0 .
\endaligned \tag13$$

The $*$-structure can be introduced as follows (\cite{7}, \cite{8},
\cite{9})
$$\chi^* (x)=\overline{\chi (S^{-1}(x^*))} . \tag14$$
It is not difficult to show that, due to the relation $S(x)^*=S^{-1}(x^*)$
valid in our $\operatorname{E(2)}$ algebra, the conjugation rule (14)
can be rephrased in Woronowicz formalism as
$$\langle \omega ,\chi^* \rangle =
-\overline{\langle \omega^*, \chi \rangle }.\tag15$$
In our case
$$\aligned
& \chi_0^*=\chi_0 +\frac{i}{\varkappa} (\chi_+ -\chi_-), \\
& \chi_{\pm}^*=-\chi_{\mp} ,\\
& {\widetilde \chi}_0^* ={\widetilde \chi}_0.
\endaligned \tag16$$
Now, it is straightforward to check that equations (13) and (16)
are equivalent to equations (11a) and (11b) provided we make
the following substitutions~:
$$\aligned
& \chi_0=e^{-\frac{P_1}{2\varkappa}}
\left( J+\frac{i}{4\varkappa} P_2 \right) ,\\
& {\widetilde \chi}_0 =-\frac{1}{4}
\left( 4\varkappa^2 \sinh^2 \left( \frac{P_1}{2\varkappa} \right)
+P_2^2 \right) ,\\
& \chi_- -\chi_+=P_2 e^{-\frac{P_1}{2\varkappa}} , \\
& I+\frac{i}{\varkappa} (\chi_+ +\chi_- )
-\frac{2}{\varkappa^2} {\widetilde \chi}_0
=e^{-\frac{P_1}{\varkappa}} .
\endaligned \tag17$$
Obviously, we should also compare the coalgebra structures. The coproduct for
$\chi$'s is defined as follows \cite{5} . Let
$\omega_i=(\omega_0,{\widetilde \omega}_0,\omega_+, \omega_-)$;
we define the functionals $f_{ij}$ by
$$\omega_i x=\sum_j (f_{ij}*x)\omega_j .\tag18$$
Then the coproduct for $\chi_i$ is defined by
$$\triangle \chi_i=\sum_j \chi_j \otimes f_{ji}
+I\otimes \chi_i .\tag19$$
We have calculated the values of $\chi_i$ and $f_{ij}$ for the generic
elements; the resulting formulae are rather involved and will not
be written out here. Then we have checked that equation (19) coincides
with equation (11c) for $\chi_0$ and $\chi_+ -\chi_-$.

\head IV.~Differential calculus on the quantum plane \endhead

Let us define the quantum plane $\varPi$ as the algebra with unity
generated by two elements $x_+$, $x_-=(x_+)^*$ subject to the following
condition
$$[x_+,x_-]=\frac{i}{\varkappa}(x_- -x_+) .\tag20$$

If we supply $\varPi$ with the coproduct, antipode and counit~:
$$\aligned
& \triangle x_{\pm}=x_{\pm} \otimes I
+I\otimes x_{\pm} ,\\
& S(x_{\pm})=-x_{\pm}, \\
& \varepsilon (x_{\pm}) =0
\endaligned \tag21$$
it becomes a quantum subgroup of deformed $\operatorname{E(2)}$.
Therefore the Woronowicz construction for $\varPi$ considered as
quantum group can be obtained by putting
\linebreak
$a\to I$, $v_{\pm}\to x_{\pm}$ in all formulae above.
As a result we obtain differential calculus determined by
the ideal generated by
$\dsize x_+^2-\frac{i}{\varkappa} x_+$,
$\dsize x_-^2 -\frac{i}{\varkappa} x_-$.
The basic rules are as follows~:
$$\aligned
& x_{\pm} \di x_{\pm}=\di x_{\pm} x_{\pm}
-\frac{i}{\varkappa} \di x_{\pm}, \\
& x_{\pm} \di x_{\mp}= -\di x_{\pm}x_{\mp}
-\frac{i}{\varkappa} \di x_{\pm}
+\di \left( x_+ x_-+\frac{i}{\varkappa}x_+ \right) .
\endaligned \tag22$$
The invariant forms read~:
$$\aligned
& \Omega_{\pm}=\di x_{\pm}, \\
& {\widetilde \Omega}_0=
\di \left( x_+ x_-+\frac{i}{\varkappa}x_+ \right)
-x_+\di x_- -x_-\di x_+ .
\endaligned \tag23$$
The external algebra reads~:
$$\aligned
& \Omega_+ \wedge \Omega_-
=-\Omega_- \wedge \Omega_+ ,\\
& {\widetilde \Omega}_0 \wedge \Omega_{\pm}
=-\Omega_{\pm} \wedge {\widetilde \Omega}_0 ,\\
& \Omega_{\alpha} \wedge \Omega_{\alpha}= 0
\qquad \alpha=+,-,0 .
\endaligned \tag24$$

Let us define the action of $\operatorname{E(2)}$ on $\varPi$ as follows~:
$$\aligned
& \varrho (I)=I\otimes I, \\
& \varrho (x_+)=a\otimes x_+ +v_+\otimes I ,\\
& \varrho (x_-)=a^*\otimes x_- + v_-\otimes I
\endaligned \tag25$$
extended by linearity and multiplicavity. Obviously this action
is a homomorphism $\varPi \to \operatorname{E(2)} \otimes \varPi$ and
obeys
$$\aligned
(\operatorname{id} \otimes \varrho )\circ \varrho
= & (\triangle \otimes \operatorname{id} )\circ \varrho, \\
(\varepsilon \otimes \operatorname{id})\circ \varrho
= & \operatorname{id}.
\endaligned \tag26$$

Let us now make the following remark. One can easily extend
the Woronowicz theory as follows.

Assume $\Cal A$ is an algebra with unity, $\Cal B$ --
quantum group and $\varrho :\Cal A\to \Cal B \otimes \Cal A$
homomorphism obeying (26). We define linear map
$$\ti \varrho :\Cal A\otimes \Cal A \to \Cal B \otimes \Cal A
\otimes \Cal A$$
as follows~: let
$$q=\sum_i x_i\otimes y_i, \qquad
\varrho (x_i)=\sum_k a_i^k\otimes x_i^k, \qquad
\varrho (y_i)=\sum_j b_i^j\otimes y_i^j,$$
then
$$\ti \varrho (q)=\sum_{i,j,k}
a_i^k b_i^j\otimes x_i^k \otimes y_i^j .\tag27$$
It is then easy to show that
$\ti \varrho :\Cal A^2\to \Cal B\otimes \Cal A^2$ and is given by~:
$$\ti \varrho \left( \sum_i x_i\otimes y_i \right)
=\sum_i \varrho(x_i) (\operatorname{id} \otimes \operatorname{D} )
\varrho (y_i) .\tag28$$
Now, we assume that $\Cal N\subset \Cal A^2$ is a sub-bimodule
such that $\ti \varrho (\Cal N) \subset \Cal B \otimes \Cal N$.
Then the differential calculus $(\Gamma, \di)$ defined by $\Cal N$
has the following property
$$\sum_k x_k\di y_k=0 \quad \Rightarrow \quad
\sum_k \varrho (x_k) (\operatorname{id} \otimes \di )
\varrho (y_k) =0 .\tag29$$
Therefore $\dsize \ti \varrho \left( \sum_k x_k \di y_k \right)
=\sum_k \varrho (x_k) (\operatorname{id} \otimes \di )
\varrho (y_k)$ is well-defined linear mapping from $\Gamma$
into $\Cal B \otimes \Gamma$.

Following the same lines as in Ref.~\cite{5} we can easily
prove the following properties of $\ti \varrho$~:
\roster
\item"{(i)}" for $x\in \Cal A$, $y\in \Gamma$
$$\aligned
& \ti \varrho (xy)=\varrho (x)\ti \varrho (y), \\
& \ti \varrho (yx)=\ti \varrho (y) \varrho (x),
\endaligned$$
\item"{(ii)}"
$$\ti \varrho \circ \di=(\operatorname{id} \otimes \di )
\circ \varrho ,$$
\item"{(iii)}"
$$\aligned
(\operatorname{id} \otimes \ti \varrho )\circ \ti \varrho
& =(\triangle \otimes \operatorname{id})\circ \ti \varrho ,\\
(\varepsilon \otimes \operatorname{id} )\circ \ti \varrho
& =\operatorname{id} .
\endaligned$$
\endroster

In our case $\Cal A=\varPi$, $\Cal B=\operatorname{E(2)}$,
$\varrho$ is given by equation (25). To check that
$\ti \varrho (\Cal N)\subset \Cal B \otimes \Cal N$ we use
the results of Ref.~\cite{5} which, together with the
property (i) above imply that it is sufficient to consider
elements $\dsize r^{-1}\left( I\otimes \left( x_{\pm}^2
-\frac{i}{\varkappa} x_{\pm} \right) \right)$.
Simple explicit calculations verify the property
under consideration
\footnote[*)]{It would probably not difficult to extend
Woronowicz theory to differential calculi on quantum
homogeneous spaces which special example we consider here.}.

The action $\ti \varrho$ is given by~:
$$\aligned
& \ti \varrho (\di x_+)= a\otimes \di x_+,  \\
& \ti \varrho (\di x_-)= a^*\otimes \di x_- ,  \\
\ti \varrho & \left(
\di \left( x_+ x_-+\frac{i}{\varkappa}x_+ \right)
-x_+ \di x_- -  x_- \di x_+ \right)  \\
& =I\otimes  \left( \di \left( x_+ x_-+\frac{i}{\varkappa}x_+ \right)
-x_+\di x_- -x_- \di x_+ \right).
\endaligned\tag30$$

Let us now describe our calculus in terms of real coordinates~:
$x_+=x_1+ix_2$, $x_-=x_1-ix_2$. We get
$$\aligned
& [x_1,x_2]=\frac{i}{\varkappa}x_2, \\
& [x_1,\di x_2]=0, \\
& [x_2, \di x_1]=-\frac{i}{\varkappa}\di x_2, \\
& [x_1,\di x_1]=\frac{1}{\varkappa^2}\varPhi, \\
& [x_2, \di x_2]=\frac{1}{\varkappa^2}\varPhi
+\frac{i}{\varkappa}\di x_1
\endaligned \tag31$$
to be compared with $3d$ calculus on $2d$ deformed Minkowski
space given by Sitarz \cite{3}.

Finally, let us construct the infitesimal transformations.
For any linear functional $\chi$ we define
$$\chi_{\varrho}=(\chi \otimes \operatorname{id}) \circ
\varrho ,\tag32a$$
$$\chi_{\ti \varrho}=(\chi \otimes \operatorname{id}) \circ
\ti \varrho .\tag32b$$
The first definition, equation (32a), coincides with the one
used by Majid and Ruegg \cite{4}, while the second one
is equivalent to Sitarz proposal, equation (23) of Ref.~\cite{3}.
Indeed, we have
$$\aligned
\chi_{\ti \varrho}(x\di y) & = (\chi \otimes \operatorname{id})
(\varrho (x)(\operatorname{id} \otimes \di )\varrho (y)) \\
& = [(\chi_{(1)}\otimes \operatorname{id})\varrho(x)]
(\operatorname{id} \otimes \di )
[(\chi_{(2)}\otimes \operatorname{id})\varrho(y)] \\
& =\chi_{(1)\varrho}(x)
(\operatorname{id}\otimes \di )
\chi_{(2)\varrho}(y) ,
\endaligned \tag33$$
where $\triangle \chi=\chi_{(1)}\otimes \chi_{(2)}$.

\enddocument